\documentclass[11pt]{article}\topmargin=-0.5cm\hfuzz=10pt  \oddsidemargin=0.3cm \textheight 220mm \textwidth=16.0cm
\usepackage{amssymb}
\newcommand{\comm}[1]{}

\def\xxxonly{\rd}
\def\xxxonly{}

\def\noxxx{\comm}
\def\Hide{\nd}

\usepackage{color}
\newcommand{\nd}[1]{\color{blue}#1\color{black}} 
\renewcommand{\nd}{\comm}
  
   \csname @addtoreset\endcsname{equation}{section}
\newtheorem{theorem}{Theorem}
\newtheorem{corollary}{Corollary}
\newtheorem{example}{Example}

\def\e{\varepsilon}

\def\defi{\stackrel{{\scriptscriptstyle \Delta}}{=}}

\def\a{\alpha}

\def\o{\omega}

\def\F{{\cal F}}
\def\w{\widehat}
\def\Ind{{\,\rm Ind\,}}
\def\Ind{{\mathbb{I}}}

\def\esssup{\mathop{\rm ess\, sup}}

\def\Re{{\rm Re\,}}

\def\R{{\bf R}}

\def\L{L}

\def\b{\beta}

\def\C{{\bf C}}

\def\ww{\widetilde}

\def\oo{\bar}

\def\L{{\cal L}}

\newcommand{\be}{\begin{equation}}
\newcommand{\ee}{\end{equation}}
\newcommand{\bd}{\begin{displaymath}}
\newcommand{\ed}{\end{displaymath}}
\newcommand{\ba}{\begin{array}{ll}}
\newcommand{\ea}{\end{array}}
\newcommand{\baa}{\begin{eqnarray}}
\newcommand{\eaa}{\end{eqnarray}}
\newcommand{\baaa}{\begin{eqnarray*}}
\newcommand{\eaaa}{\end{eqnarray*}}

\def\oo{\bar}

\def\a{\alpha}

\def\HH{\mathbb{H}}
 
\def\ew{\left(i\o\right)}

\date{ }
\title{On  incompleteness of polynomials in some weighted spaces on half line }
\author{Nikolai Dokuchaev}
 
\begin{document}
\def\brea{}
\def\breakk{}
\def\break{}
\def\breakm{\nonumber\\  }
\def\breacm{}
\def\BRR{}
\date{Submitted: June 1 2020. Revised: November 4 2020}
\maketitle
\noxxx{ \let\thefootnote\relax\footnote{The author is with
Zhejiang University/University of Illinois at Urbana-Champaign Institute,  Zhejiang University, Haining, China.}}
\begin{abstract} 
The paper studies  completeness of the polynomials in 
weighted $L_p$-spaces on half line. It is shown  that the completeness of polynomials does not hold for a wide class of weights, including the weights  $\exp(- r t^q)$ with $r>0$ and $q\in (0,1)$. 
\end{abstract}
\section{Introduction}
The theory of  approximation of function by polynomials is well developed; however, the exact characterisation of the weights for which polynomials 
are complete is unknown for weighted $L_2$-spaces of functions  on infinite intervals. 
The related questions have been studied intensively;   see, e.g., \cite{Bakan,BC,Higgins,Natan,Sc},
and the literature therein.

\comm(
For example, it is known 
that the moment problem is indeterminate in the weighted $L_2$-space of functions defined on the entire line $\R$ 
with the weight $\rho$ such that the following Krein condition holds:\baa
 \int_{-\infty}^\infty\frac{ \log \rho(\o)}{1+\o^2}d\o>-\infty;
\label{K0}\eaa
see, e.g., Theorem 4.14 \cite{Sc}.  This implies that polynomials 
are  not complete in this weighted space; see Theorems 6.10 and 7.7 in \cite{Sc}.

It is also known that  the Stieltjes moment problem is indeterminate   in the weighted $L_2$-space of functions defined on $[0,+\infty)$ 
with the weight $\rho$ such that the following Krein condition holds:\baa
 \int_{0}^\infty\frac{ \log \rho(\o^2)}{1+\o^2}d\o =\int_{0}^\infty\frac{ \log \rho(\o)}{(1+\o)\sqrt{\o}}d\o>-\infty;
\label{K1}\eaa
see, e.g., Theorem 4.17 \cite{Sc} or Theorem 1.2 \cite{p2}. 
In Theorem 4.17 \cite{Sc}, this condition is stated as a sufficient condition of indeterminacy; in Theorem 1.2 \cite{p2}, this condition is stated as a necessary and sufficient  condition of indeterminacy.
Theorem 25 in \cite{Merg} gives necessary and sufficient conditions of completeness of polynomials  on the half line
with respect to the norm $\|x\|=\sup_\o |x(\o)|\rho(\o)$ that is not equivalent to the weighted  $L_p$-norms considered in the present paper. 
Some developments  and historical notes can be found \cite{Bakan,Merg,p1,p2,Sc}; see also references therein.

The present paper focuses on the problem of  completeness of polynomials on half line in a framework that bypasses  the  moment problem and indeterminacy. It is shown directly  that  polynomials are incomplete on a weighted $L_p$-spaces on half line under a weaker condition than (\ref{K1}). In particular, this condition allows weights  $\exp(- r \o^q)$ with $r>0$ and $q\in [1/2,1)$ that are excluded by condition (\ref{K1}).

\section{The result}\label{secSetting}
Let $\R^+\defi [0,+\infty)$, and let $\cal R^+$  be the set of measurable functions $\rho:\R^+\to \R^+$  such that $\int_0^{+\infty}\o^k \rho(\o)d\o<+\infty$,
for all $k=0,1,2,....$, and
\baa
\int_0^\infty\frac{ \log \rho(\o)}{1+\o^2}d\o>-\infty.
\label{K2}\eaa
\begin{example} It can be verified directly that $\cal R^+$ includes 
 $\rho(\o)=e^{-(\log |\o|)^2}$,
$\rho(\o)=e^{-r|\o|^q}$, and $\rho(\o)=e^{-r|\o|^q/|\log |\o||^p}$ for $r>0$, $q\in (0,1]$, $p\ge 2$.
On the other hand, the function  $\rho(\o)=e^{-|\o|/|\log |\o||}$ is excluded.
\end{example}

For a  interval $I\subset \R$, for a measurable function $\varrho:I\to\R^+$,
  and for $p\ge 1$, let  $L_{p,\varrho}(I)$ be the Banach space of complex valued functions $u:\R\to\C$ with the norm 
 \baaa \|u\|_{L_{p,\varrho}(I)}=\Bigl(\int_{I}\varrho(\o)|u(\omega)|^p d\o\Bigr)^{1/p}.
 \eaaa
 
 \begin{theorem}\label{ThM}  For  any $\rho\in\cal R^+$,  $T>0$, $r>0$, and $q\in (0,1)$,
 the function $e^{i\o T}$  cannot be approximated 
 by polynomials in the space $L_{1,\rho}(\R^+)$. 
 \end{theorem}
 \begin{corollary}\label{corr1} For any  $\rho\in\cal R^+$, $r>0$ and $q\in (0,1)$,
 the set of polynomials is  incomplete  in $L_{p,\rho}(D)$ for all $p\ge 1$.
 \end{corollary}
 
 The integrand in condition (\ref{K2}) is the same as the one in the  Krein condition  for the  entire real line (\ref{K0}).
  This condition is less restrictive than condition (\ref{K1}). 
For example, 
for $\rho(\o)=\exp(- r \o^q)$ with $r>0$ and $q\in [1/2,1)$,
condition (\ref{K2}) holds but condition (\ref{K1}) does not hold. 
\Hide{see Example 4.19 in \cite{Sc}.} Furthermore, condition (\ref{K2})  is different from  the condition from Theorem 25
\cite{Merg}, where completeness was considered with respect to the norm $\|x\|=\sup_{\o}\rho(\o)|x(\o)|$, which is not equivalent
to the norms in the present paper.
\section{Proofs} \xxxonly{
 \par
The proof of  Theorem \ref{ThM} below is based on the approach \cite{D20}
developed for  analysis of predictability of processes with fast decaying Fourier transform.   }
\subsection{Some background notations}
 Let $\R^-\defi\{\o\in\R: \o<0\}$,   $\C^+\defi\{z\in\C:\
\Re z> 0\}$, $\C^-\defi\{z\in\C:\
\Re z< 0\}$, $i=\sqrt{-1}$.
\par
For $p\in[1,+\infty]$ and intervals $I\subset \R$, we denote by $L_p(I)$   the usual $L_p$-spaces
of functions  $x:I\to \C$.

\par
For $x\in  L_p(\R)$, $p=1,2$, we denote by $X=\F x$ the function
defined on $i\R$ as the Fourier transform of $x$; $$X(i\o)=(\F
x)(i\o)= \int_{-\infty}^{\infty}e^{-i\o t}x(t)dt,\quad \o\in\R.$$ If
$x\in L_2(\R)$, then $X$ is defined as an element of $L_2(i\R)$, i.e.,  $X(i\cdot)\in L_2(\R)$.
\par
For $x(\cdot)\in L_p(\R)$, $p=1,2$, such that $x(t)=0$ for $t<0$, we denote by
$\L x$  the Laplace transform \baa\label{Up} X(z)=(\L
x)(z)\defi\int_{0}^{\infty}e^{-z t}x(t)dt, \quad z\in\C^+. \eaa
In this case, $X|_{i\R}=\F x$.
\par
Let  $\HH^r$ be the Hardy space of holomorphic on $\C^+$ complex valued functions
$h(z)$ with finite norm
$\|h\|_{\HH^r}=\sup_{s>0}\|h(s+i\cdot)\|_{L_r(\R)}$, $r\in[1,+\infty]$; see, e.g., \cite{Du}. 

Similarly to ${\cal R}^+$, we denote by $\cal R$  the set of measurable functions $\rho:\R\to [0,+\infty)$  such that, for $k=0,1,2,....$, 
$\int_{-\infty}^{+\infty}|\o|^k \rho(\o)d\o<+\infty$, and 
\baaa
\int_{-\infty}^\infty\frac{ \log \rho(\o) }{1+\o^2}d\o>-\infty.
\eaaa

\subsection{A supposition}

We assume below that we are given $T>0$  and $\rho\in{\cal R}^+$.

Suppose that the theorem statement is incorrect. Then there exists 
sequence of  polynomials  $\{\ww\psi_d(\o)\}_{d=1}^\infty$ in $\o\in \R$ of order $d$ such that 
\baa
&&\e_d\defi \|e^{iT\cdot}-\ww\psi_d(\cdot)\|_{L_{1,\rho}(\R^+)}\breakk =\int_{0}^\infty |e^{iT\o}-\ww\psi_d(\o)|  \rho(\o)d\o\to 0\quad \breakk \hbox{as}\quad d\to +\infty. \quad
\label{lim1}
\eaa



\subsection{Some featured functions}
\subsubsection{Weights $\rho_d$}
We presume that the given $\rho\in{\cal R}^+$ is extended to  $\rho\in{\cal R}$; a possible choice is  such that $\rho$ is an even function. 

Let
\baaa
&&L_d\defi \|e^{iT\cdot}-\ww\psi_d(\cdot)\|_{L_{1,\rho}(\R^-)}^2\breakk =\int^{0}_{-\infty} |e^{iT\o}-\ww\psi_d(\o)|  \rho(\o)d\o. 
\label{Ld}
\eaaa

Consider a sequence   of functions $\{\rho_d(\o)\}_{d=1}^\infty$ defined as 
\baaa
\rho_d(\o)=\Ind_{\{\o<0\}} \rho(\o)\frac{\e_d}{L_d}+\Ind_{\{\o\ge 0\}}\rho(\o).
\eaaa
Since $\e_d>0$ for all $d$,  we have that $\{\rho_d(\o)\}_{d=1}^\infty\subset {\cal R}$.
\subsubsection*{Functions $X_d$}
Let us construct a    sequence   of 
functions $\{X_d(\o)\}_{d=1}^\infty\subset  \HH^2$ such that 
$|X_d( i \o)|=\mu_d(\o)$, where 
\baaa
&&\mu_d(\o)\defi\rho_d(\o)/(1+\o^2).
\eaaa
  Existence of such  functions
follows from Theorems  11.6 and 11.7 from \cite{Du}. For example, one 
can select  
\Hide{\baaa
In\ Duren\ \
Y(p)=\exp\left[\frac{1}{\pi i}\int_{-\infty}^\infty 
\frac{(1+s p)\log \mu(s)}{(s-p)(1+s^2)}ds\right],\quad z\in \C.
\eaaa
$X(z)=Y(-iz)$}
\baaa
X_d(z)=\exp\left[\frac{1}{\pi i}\int_{-\infty}^\infty 
\frac{(1-isz)\log \mu_d(s)}{(s+iz)(1+s^2)}ds\right],\quad z\in \C.
\eaaa
We used here equation (11) from Theorem 11.6 \cite{Du} stated
for the Hardy spaces on the upper complex half-plane; in the present paper, it is adjusted to the Hardy spaces $\HH^p$ on the right  half-plain. 
In particular, we have that 
\baaa
\int_{-\infty}^\infty\frac{ |\log|X_d(i\o)|| }{1+\o^2}d\o=
\int_{-\infty}^\infty\frac{ |\log \rho_d(\o)|}{1+\o^2}d\o\brea+
\int_{-\infty}^\infty\frac{ |\log ((1+\o^2)^{-1})| }{1+\o^2}d\o<+\infty.
\eaaa

Let $x_d\defi \F^{-1}X_d|_{i\R}$; clearly, $x_d(t)=0$ for $t<0$.  

Clearly,
\baa
\sup_d\|X_d(i\cdot)\|_{L_2(\R)}<+\infty, \quad  \sup_d\|x_d\|_{L_2(\R)}<+\infty.
\label{XX}
\eaa

\subsubsection*{Convolution kernels $h_d$ and functions $y_d$}

Let us construct  functions $h_d:\R\to\C$
such that $h_d(t)=0$ for $t\notin [-T,0]$ , $h_d\in
C^\infty(\R)$, and 
\baa
\inf_d\left|\int_{0}^{T} h_d(-t)x_d(t)dt\right|>0. 
\label{y0}\eaa

First,  let us define
 $\varkappa_\e (t)\defi\e^{-1}\varkappa_1(t/\e )$,
where $\varkappa_1(t)\defi\exp(t^2(1-t^2)^{-1})$ is the
so-called Sobolev kernel. Let \baaa
 g_d(t)\defi\|x_d|_{[0,T]}\|_{L_2(0,T)}^{-2}x_d(-t).
 \eaaa
Finally, let $h_{d,\e}$ be defined as the convolution  
\baaa
h_{d,\e}(t)=\int_{-\infty}^{\infty}\varkappa_\e (t-s)\Ind_{[-T+\e,-\e]}(s)g_d(s)ds,\quad \e>0.\label{filter}
\eaaa
In this case, we have that $h_{d,\e}\to g_d$ in $L_2(\R)$ as $\e\to 0$.

It follows from (\ref{XX}) that
\baa
\sup_d\|h_d(i\cdot)\|_{L_p(\R)}<+\infty, \qquad  p=1,2.
\label{hh}
\eaa

Furthermore, let
\baa
y_d(t)\defi \int_{t}^{t+T} h_d(t-s)x_d(s)ds.
\eaa

 Let $\xi_{d,\e}(t)\defi h_{d,\e}(-t)-g_d(-t)$.
We have that $h_d(-t)=h_{d,\oo\e}(-t)=g_{d}(-t)+\xi_{d,\oo\e}(t)$,
\baaa
&&y_d(0)=\int_{0}^{T} h_d(-t)x_d(t)ds\\&&=\|x_d|_{[0,T]}\|_{L_2(0,T)}^{-2}
\int_{0}^{T} |x_d(t)|^2dt\breakk+ \|x_d|_{[0,T]}\|_{L_2(0,T)}^{-2}
\int_{0}^{T} \xi_{d,\e}(t)x_d(t) dt \\
&&=1+\|x_d|_{[0,T]}\|_{L_2(0,T)}^{-2}\int_{0}^{T} \xi_{d,\e}(t) x_d(t) dt.
\eaaa
Furthermore,
\baaa
&&|y_d(0)-1|\breakk \le \|x_d|_{[0,T]}\|_{L_2(0,T)}^{-2}\|x_d|_{[0,T]}\|_{L_2(0,T)}\|\xi_{d,\e}|_{[0,T]}\|_{L_2(0,T)}
\breakk=\|x_d|_{[0,T]}\|_{L_2(0,T)}^{-1}\|\xi_{d,\e}|_{[0,T]}\|_{L_2(0,T)}.
\eaaa
We have that the sequence $\{X_d\}$ has a limit in $L_2(\R)$, 
therefore, the  sequence $\{x_d\}$ has a limit in $L_2(\R)$, and 
\baaa
\|x_d|_{[0,T]}\|_{L_2(0,T)}^{-1}\|\xi_{d,\e}|_{[0,T]}\|_{L_2(0,T)}\le C
\eaaa
 for all $d,\e$ for some $C>0$.

By the property of  kernels $\varkappa_\e$, for any $d$,
\baaa
\|x_d|_{[0,T]}\|_{L_2(0,T)}^{-1}\|\xi_{d,\e}|_{[0,T]}\|_{L_2(0,T)}\to 0 \quad \hbox{as}\quad \e\to 0.
\eaaa
Hence, for each $d$,  for some  choice of $\oo\e=\oo\e(d)>0$, the kernels selected as  $h_d=h_{d,\oo \e}$  are such that 
\baa
\inf_d|y_d(0)|> 0.
\label{yneq0}
\eaa
It follows that (\ref{y0}) holds for this choice of $h_d$.


\subsection{Proof of Theorem \ref{ThM}: final steps}
Let  $H_d\defi \L h_d$ and $Q_d=\L q_d$,  where $q_d(t)\equiv h_d(t-T)$. 
Clearly,  $H_d\ew\equiv Q_d\ew e^{i\o T}$, 
 $q_d(t)=0$ for $t<0$, $q_d\in C^{\infty}(\R)$, and  $Q_d\in \HH^2\cap\HH^\infty$. 

By the choice of $\rho_d$, it follows that
\baaa
\|\ww\psi_d(\cdot)-e^{i\cdot T}\|_{L_{1,\rho_d}(\R)} =2\|\ww\psi_d(i\cdot)-e^{i\cdot T}\|_{L_{1,\rho_d}(\R^+)}=2\e_d\to 0\quad
\breakk \hbox{as}\quad d\to +\infty.\label{lim20}
\eaaa

For our purposes, it would be more convenient to use  a
sequence of  polynomials  $\{\psi_d(z)\}_{d=1}^\infty$  of order $d$ such that 
 \baa
\|e^{iT\cdot}-\psi_d(i\cdot)\|_{L_{2,\rho_d}(\R)}^2\breakk
=2\|\ww\psi_d(\cdot)-e^{i\cdot T}\|_{L_{1,\rho_d}(\R^+)} \rho_d(\o)d\o\to 0\quad
\breakk \hbox{as}\quad d\to +\infty.\label{lim2}\eaa 
The coefficients $a_k$ of the polynomials $\psi_d(z)=\sum_{k=0}^d a_k z^k$ can be constructed by adjustment the signs of the coefficients for the polynomials
$\ww\psi_d(\o)=\sum_{k=0}^d \ww a_k \o^k$ such that $\psi_d(i\o)\equiv \ww\psi_d(\o)$, i.e., $\ww a_k=a_k i^k$
and $a_k=\ww a_k i^{-k}$. 
 In this case, statements (\ref{lim1}) and (\ref{lim2}) are equivalent.

For $d=1,2,....$, set 
\baa &&\w H_d(z)
\defi e^{-Tz} \psi_d (z) H_d(z)=\psi_d(z)Q_d(z),\quad\breakk\w h_d\defi \F^{-1}\w H_d|_{i\R}.
\label{Kk}
\label{FD} \eaa

Since  $q_d(t)=0$ for $t<0$ and $q_d\in C^{\infty}(\R)$, we have that $z^n Q_d(z)\in \HH^2\cap\HH^\infty$ for any integer $n\ge 0$.
It follows  that \baaa
\w H\in \HH^2\cap\HH^\infty.
\label{1}\eaaa

Let   $\w Y_d(i\o)\defi \w H_d(i\o)X_d(i\o)$ and 
 $\w y_d\defi \F^{-1}\w Y_d(i\o)$.  It follows that
\baa
\w y_d(t)=\int^t_{-\infty}\w h_d(t-s)x_d(s)ds.\label{TD}
\eaa

It can be noted that, by the definitions, 
$
\w h_d(t)\defi \sum_{k=0}^d a_{dk}\frac{d^kh}{dt^k}(t+T).
$
\par
We have that
\baa
&&\|\w y_{d}-y_d\|_{L_\infty(\R)}\breakk \le \frac{1}{2\pi} \|(\w H_d(i\cdot) -H_d(i\cdot))X(i\cdot)\|_{L_1(\R)}.
\label{est1}\eaa
Furthermore, 
\baa
&&\|(\w H_d(i\cdot) -H_d(i\cdot))X_d(i\cdot)\|_{L_1(\R)}\breakk=\int_{-\infty}^\infty \Bigl|(e^{- i\o T}\psi_d(i\o)-1) 
e^{ i\o T} Q_d(i\o)X_d(i\o)\Bigr|d\o
\nonumber \\
&&= \int_{-\infty}^\infty \rho_d(\o)\Bigl| (e^{- i\o T}\psi_d(i\o)-1) \BRR
\rho_d(\o)^{-1}e^{ i\o T} Q_d(i\o)X_d(i\o)\Bigr|d\o 
\nonumber\\
&&=\int_{-\infty}^\infty \rho_d(\o) \Bigl|(\psi_d(i\o) -e^{i\o T}) \BRR
\rho_d(\o)^{-1}  e^{ i\o T} Q_d(i\o)X_d(i\o)\Bigr|d\o \breakk \le \a_d  \b_d. \label{est}
\eaa
Here
\baaa
&&\a_d\defi \int_{-\infty}^\infty \rho_d(\o)|\psi_d(i\o)-e^{i\o T}|d\o\breakk=\|\psi_d(i\cdot)-e^{i\cdot T}\|_{L_{1,\rho_d}(\R^-)}+\|\psi_d(i\cdot)-e^{i\cdot T}\|_{L_{1,\rho_d}(\R^+)},
\eaaa
and \baaa
&& \b_d\defi \esssup_\o \rho_d(\o)^{-1}|e^{ i\o T} Q_d(i\o)X_d(i\o)|.
\eaaa
Remind that, by the choice of $\rho_d$, it follows that \baaa
\|\psi_d(i\cdot)-e^{i\cdot T}\|_{L_{1,\rho_d}(\R^-)} =\|\psi_d(i\cdot)-e^{i\cdot T}\|_{L_{1,\rho_d}(\R^+)}=\e_d.
\eaaa
It gives  that
\baa
 \a_d\to 0\quad \hbox{ as}\quad d\to +\infty.
\label{a} 
\eaa
\par
 It follows from estimate (\ref{hh}) with $p=1$ that
\baaa
\sup_{\o,d} |Q_d(i\o)|<+\infty.
\label{QQ}
\eaaa
Hence, by the choice of $X_d$, we obtain that
\baa
\sup_d |\b_d|\le \esssup_{\o,d} |Q_d(i\o)| \rho_d(\o)^{-1} |X_d(i\o)|\brea
=\esssup_{\o,d} |Q_d(i\o)| (1+\o^2)^{-1}<+\infty.
\label{b}\eaa Then estimates (\ref{est1})--(\ref{b}) imply 
that 
\baa
&&\sup_t|\w y_{d}(t)-y_d(t)|\to 0\quad \hbox{ as}\quad d\to +\infty.
\label{est2}\eaa

On the other hand,  since  $x_d(t)$ and $\w h_d(t)$  are both vanishing for $t<0$, we have that
\baa
\w y_d(t)=0, \quad t\le 0.
\label{est3}\eaa
Hence (\ref{est2}) cannot hold   simultaneously with (\ref{yneq0}) and (\ref{est3}).  
This means that the supposition that the theorem statement is incorrect leads to a contradiction.  
This completes the proof of Theorem \ref{ThM}. $\Box$.

{\em The proof of Corollary \ref{corr1}} follows from the fact that $L_{p,\rho}(\R^+)\subseteq
L_{1,\rho}(\R^+)$ for any $p\ge 1$, and that this embedding is continuous.

\subsection*{Acknowledgement}
The author thanks Prof.
Grzegorz \'Swiderski  for his valuable comments and advices that helped to improve the paper.

 
\end{document}